\documentclass{IEEEtran4PSCC}

\usepackage{cite}

\ifCLASSINFOpdf
   \usepackage[pdftex]{graphicx}
\else
   \usepackage[dvips]{graphicx}
\fi
\usepackage[cmex10]{amsmath}
\usepackage{amsfonts}
\usepackage{xcolor}
\usepackage{tikz}
\usetikzlibrary{calc}
\usetikzlibrary{shapes}
\usetikzlibrary{arrows}
\usetikzlibrary{positioning}
\usepackage{circuitikz}
\usepackage{enumerate}
\usepackage{enumitem}
\usepackage{algpseudocode}
\usepackage{algorithm}
\usepackage{empheq}
\usepackage{tikz,pgfplots}

\makeatletter
\let\old@ps@headings\ps@headings
\let\old@ps@IEEEtitlepagestyle\ps@IEEEtitlepagestyle
\def\psccfooter#1{%
    \def\ps@headings{%
        \old@ps@headings%
        \def\@oddfoot{\strut\hfill#1\hfill\strut}%
        \def\@evenfoot{\strut\hfill#1\hfill\strut}%
    }%
    \def\ps@IEEEtitlepagestyle{%
        \old@ps@IEEEtitlepagestyle%
        \def\@oddfoot{\strut\hfill#1\hfill\strut}%
        \def\@evenfoot{\strut\hfill#1\hfill\strut}%
    }%
    \ps@headings%
}
\makeatother

\begin{document}
%
\title{Tight Big-Ms for Optimal Transmission Switching}

\author{
\IEEEauthorblockN{Salvador Pineda, Juan Miguel Morales, Álvaro Porras, Concepción Domínguez}
\IEEEauthorblockA{OASYS research group, University of Málaga, Spain \\
\{spineda, juan.morales, alvaroporras, concepcion.dominguez\}@uma.es}
}


\maketitle

\begin{abstract}
This paper addresses the Optimal Transmission Switching (OTS) problem in electricity networks, which aims to find an optimal power grid topology that minimizes system operation costs while satisfying physical and operational constraints. Existing methods typically convert the OTS problem into a Mixed-Integer Linear Program (MILP) using big-M constants. However, the computational performance of these approaches relies significantly on the tightness of these big-Ms. In this paper, we propose an iterative tightening strategy to strengthen the big-Ms by efficiently solving a series of bounding problems that account for the economics of the OTS objective function through an upper-bound on the generating cost. We also discuss how the performance of the proposed tightening strategy is enhanced if reduced line capacities are considered. Using the 118-bus test system we demonstrate that the proposed methodology outperforms existing approaches, offering tighter bounds and significantly reducing the computational burden of the OTS problem. 
\end{abstract}


\thanksto{\noindent This work was supported in part by the European Research Council (ERC) under the EU Horizon 2020 research and innovation program (grant agreement No. 755705), in part by the Spanish Ministry of Science and Innovation (AEI/10.13039/501100011033) through project PID2020-115460GB-I00. \'A. Porras is also financially supported by the Spanish Ministry of Science, Innovation and Universities through the University Teacher Training Program with fellowship number FPU19/03053. Finally, the authors thankfully acknowledge the computer resources, technical expertise, and assistance provided by the SCBI (Supercomputing and Bioinformatics) center of the University of M\'alaga.}

\section{Introduction} \label{sec:intro}

Traditionally, transmission lines in electricity networks have been regarded as infrastructure devices that cannot be controlled, except during instances of outages or maintenance. More recently, the possibility of flexibly exploiting the topological configuration of the grid was first suggested in~\cite{o2005dispatchable} and later formalized in~\cite{fisher2008optimal} into what we know today as the \emph{Optimal Transmission Switching} (OTS) problem. The optimal transmission switching problem refers to the task of determining the most efficient configuration of transmission lines in a power system to achieve certain objectives. It involves deciding which transmission lines should be open or closed to optimize system performance in terms of factors such as minimizing transmission losses, voltage deviations, or congestion. Even if the power flow equations are simplified using the well-known \emph{direct current} (DC) linear approximation of the power flow equations, the resulting formulation of the OTS problem, known as DC-OTS, takes the form of a mixed-integer program, which has been proven to be NP-hard for general network classes~\cite{kocuk2016cycle, fattahi2019bound}. 

Up until now, the resolution of the DC-OTS has been addressed using two different methodological approaches. These approaches can be categorized as \textit{exact} methods and \textit{heuristics}. The exact methods utilize techniques derived from mixed-integer programming, such as bounding  and generating valid cuts. These methods aim to solve the DC-OTS with (certified) global optimality, ensuring the best possible solution. On the other hand, heuristics aim to rapidly identify good solutions for the problem, potentially sacrificing optimality or even  suggesting infeasible grid configurations.

Several heuristic methods have been proposed in the technical literature to reduce the computational time in solving the OTS problem. Some of these methods focus on decreasing the number of lines that can be switched off \cite{liu2012heuristic, barrows2012computationally, flores2020alternative}. Other approaches maintain the original set of switchable lines but determine their on/off status using greedy algorithms \cite{fuller2012fast, crozier2022feasible}. In contrast, the authors of \cite{hinneck2022optimal} propose a parallel approach where heuristics generate promising candidate solutions to expedite traditional MIP algorithms in solving the OTS problem. Furthermore, certain data-based heuristic methods utilize information from past OTS problems to improve efficiency. For example, the authors of \cite{johnson2020knearest, pineda2023learning} employ a $K$-nearest neighbor strategy to significantly reduce the search space of the integer solution for the DC-OTS problem. Similarly, references \cite{yang2019line, han2022learning} present more advanced methodologies involving neural networks to learn the optimal status of switchable lines.

Within the exact methods, notable contributions can be found in references \cite{kocuk2016cycle}, \cite{fattahi2019bound}, and \cite{dey2022node}. In particular, the authors of \cite{kocuk2016cycle} present a cycle-based formulation for the DC-OTS problem, which yields a mixed-integer linear program. They also introduce sets of strong valid inequalities for a relaxed version of their formulation that can be efficiently separated. In \cite{fattahi2019bound}, the authors focus on the mixed-integer linear formulation of the DC-OTS, utilizing a big-M approach to handle the disjunctive relationship between the power flow of switchable lines and the voltage angle differences. They prove the NP-hardness of determining valid big-Ms and propose a methodology for setting the appropriate values. Lastly, the authors of \cite{dey2022node} develop a family of cutting planes specifically tailored for the DC-OTS problem.

This paper introduces a new exact methodology to address the DC-OTS problem, making significant contributions to the existing state-of-the-art. Our approach determines suitable values for the big-M constants used in the mixed-integer reformulation of the DC-OTS by solving the so-called \textit{bounding problems}. To obtain tighter big-M values, we impose an upper-bound on the generating cost in these bounding problems. We also investigate the synergistic effect between the big-M tightening and the reduction of the line capacities. The performance of our methodology is then compared to the approach proposed in \cite{fattahi2019bound} to determine big-M constants for the 118-bus test system. In summary, the key contributions of this work can be summarized as follows:
\begin{itemize}
\item[-] We propose a set of bounding problems to efficiently compute tight big-M values to be used in the optimal transmission switching problem. Besides, we enhance the performance of the bounding problems imposing a valid upper-bound on the total generating cost.
\item[-] We extend the use of the bounding problem to also compute the maximum feasible power flow through the transmission lines. We demonstrate that using these reduced capacities decrease even further the big-M values.
\item[-] We use a 118-bus test system to prove that the proposed bound tightening methodology clearly outperforms state-of-the-art approaches and general purpose methods of optimization solvers in terms of the computational burden required to solve the OTS problem. 
\end{itemize}

The remainder of this paper is structured as follows. Section \ref{sec:ots} introduces the original formulation of the DC-OTS problem, its reformulation as a mixed-integer linear program, and the existing methodologies to compute the required big-M constants. The proposed cost-driven bound tightening approach is presented in Section \ref{sec:methodology}, which concludes with the comparison procedure used to assess its performance. Section \ref{sec:case_study} discusses the computational results obtained for the 118-bus test system. Finally, conclusions are drawn in Section \ref{sec:conclusions}.

\section{Optimal Transmission Switching} \label{sec:ots}
In this section we introduce the standard and well-known formulation of the \emph{Direct Current Optimal Transmission Switching} problem (DC-OTS) which has been considered in \cite{fisher2008optimal, kocuk2016cycle, fattahi2019bound, liu2012heuristic, barrows2012computationally, fuller2012fast, hinneck2022optimal, johnson2020knearest, pineda2023learning, yang2019line, dey2022node, hedman2012flexible, moulin2010transmission} among others. Readers interested in exploring further details regarding the additional complexity of the AC-OTS problem can refer to \cite{flores2020alternative, crozier2022feasible, numan2023role}.

Consider a power network consisting of a set of nodes $\mathcal{N}$ and transmission lines $\mathcal{L}$. For simplicity, we assume that there is one generator and one power load per node $n \in \mathcal{N}$. Let $p_n$ and $d_n$ denote the power dispatch of the generator and the power consumed by the power load, respectively. Each generator is characterized by a minimum and maximum power output, $\underline{p}_n$ and $\overline{p}_n$, and a marginal production cost $c_n$. We represent the power flow through the line $(n,m) \in \mathcal{L}$ connecting nodes $n$ and $m$ by $f_{nm}$. As customary, $f_{nm}>0$ represents a power flow from node $n$ to node $m$, and $f_{nm}<0$ a power flow in the opposite direction. The maximum power flow from node $n$ to node $m$ is denoted by $F_{nm}$, and the power flow in the opposite direction is limited by $F_{mn}$. The maximum flow through a line is constrained by thermodynamics limitations and therefore, \emph{physical} line capacities are always symmetric, i.e., $F_{nm}=F_{mn}$. However, depending on the location of generators and loads in a network, the maximum \emph{feasible} flows through a line may be different in each direction, which is the reason why we consider the more general case of asymmetric line capacities. 
Besides, the set of transmission lines that can be switched on/off is denoted by $\mathcal{L}_S \subseteq \mathcal{L}$. If the line $(n,m) \in \mathcal{L}_S$, the binary variable $x_{nm}$ determines its status, being equal to 1 if the line is fully operational, and 0 when disconnected. Using the DC approximation of the network equations, the flow $f_{nm}$ through an operational line is given by the product of the susceptance of the line, $b_{nm}$, and the difference of the voltage angles at nodes $n$ and $m$, i.e., $\theta_n-\theta_m$. We use bold symbols to define the vectors of variables $\mathbf{p}=[p_n, n \in \mathcal{N}]$, $\boldsymbol{\theta}=[\theta_n, n \in \mathcal{N}]$, $\mathbf{f}=[f_{nm}, (n,m) \in \mathcal{L}]$, and $\mathbf{x}=[x_{nm}, (n,m) \in \mathcal{L_S}]$. With this notation in place, the DC-OTS problem can be formulated as follows:
\begin{subequations}\label{eq:OTS_NP}
\begin{align}
& \min_{p_n,f_{nm},\theta_n,x_{nm}}  \quad \sum_{n} c_{n} \, p_{n} \label{eq:OTS_NP_obj}\\
\text{s.t.} & \quad f_{nm} = x_{nm}b_{nm}(\theta_n-\theta_m), \quad \forall (n,m) \in \mathcal{L}_S \label{eq:OTS_NP_Flow_S}\\
& \quad f_{nm} = b_{nm}(\theta_n-\theta_m), \quad \forall (n,m) \in \mathcal{L} \setminus \mathcal{L}_S\label{eq:OTS_NP_Flow_NS}\\
&\quad \hspace{-5mm} \sum_{m:(n,m)\in\mathcal{L}}f_{nm} - \hspace{-5mm} \sum_{m:(m,n)\in\mathcal{L}}f_{mn} = p_n - d_n, \quad \forall n \in \mathcal{N} \label{eq:OTS_NP_PB}\\
& \quad \underline{p}_n \leq p_n \leq \overline{p}_n, \quad \forall n \in \mathcal{N} \label{eq:OTS_NP_Plimits}\\
& \quad -x_{nm}F_{mn}\leq f_{nm} \leq x_{nm}F_{nm}, \quad \forall (n,m) \in \mathcal{L}_S \label{eq:OTS_NP_Flow_limit_S}\\
& \quad -F_{mn}\leq f_{nm} \leq F_{nm}, \quad \forall (n,m) \in \mathcal{L} \setminus \mathcal{L}_S\label{eq:OTS_NP_Flow_limit_NS}\\
& \quad \theta_1 = 0 \label{eq:OTS_NP_slack}\\
& \quad x_{nm} \in \{0,1\}, \quad \forall (n,m) \in \mathcal{L}_S\label{eq:OTS_NP_binary}
\end{align}
\end{subequations}

The objective function~\eqref{eq:OTS_NP_obj} minimizes the total electricity generation cost. The power flow through transmission lines is defined in~\eqref{eq:OTS_NP_Flow_S} and~\eqref{eq:OTS_NP_Flow_NS}. In the case of a switchable line, constraint \eqref{eq:OTS_NP_Flow_S} includes the binary variable $x_{nm}$ to enforce this  relationship only when the line is in service. Naturally, $x_{nm} = 0$ implies that $f_{nm} = 0$. The nodal power balance equation is ensured by~\eqref{eq:OTS_NP_PB}, while constraints \eqref{eq:OTS_NP_Plimits} impose that the power output of generating units must lie within the interval $[\underline{p}_n,\overline{p}_n]$. Constraints~\eqref{eq:OTS_NP_Flow_limit_S} and~\eqref{eq:OTS_NP_Flow_limit_NS} limit the maximum power flow through switchable and non-switchable lines, respectively. Equation~\eqref{eq:OTS_NP_slack} arbitrarily sets one of the voltage angles to zero, while the binary character of variables $x_{nm}$ is imposed by constraint~\eqref{eq:OTS_NP_binary}.  

Problem~\eqref{eq:OTS_NP} is a mixed-integer nonlinear programming problem due to the product $x_{nm}(\theta_n-\theta_m)$ in~\eqref{eq:OTS_NP_Flow_S}. Even when the power network includes a connected subgraph of non-switchable lines, this problem has been proven to be NP-hard \cite{fattahi2019bound}. However, constraint~\eqref{eq:OTS_NP_Flow_S} can be linearized by introducing a pair of large enough constants $M_{nm}$, $M_{mn}$ per switchable line \cite{hedman2012flexible}. By doing so, equation~\eqref{eq:OTS_NP_Flow_S} can be replaced by the two following inequalities:
\begin{subequations}\label{eq:linear}
\begin{align}
& f_{nm} \geq b_{nm}(\theta_n-\theta_m)-M_{nm}(1-x_{nm}) \\
& f_{nm} \leq b_{nm}(\theta_n-\theta_m)+M_{mn}(1-x_{nm}) 
\end{align}
\end{subequations}
where the large constants $M_{nm}, M_{mn}$ are guaranteed to be upper bounds of $b_{nm}(\theta_n-\theta_m)$ and $b_{nm}(\theta_m-\theta_n)$, respectively, when the line $(n,m)$ is disconnected ($x_{nm} = 0$). Under that assumption, the DC-OTS is reformulated as the following mixed-integer linear programming problem
\begin{subequations} 
\begin{align}
 \min_{p_n,f_{nm},\theta_n,x_{nm}}  & \quad  \sum_{n} c_{n} \, p_{n} \label{eq:ots_mip_of}\\
\text{s.t.}  & \quad \eqref{eq:OTS_NP_Flow_NS}-\eqref{eq:OTS_NP_binary}, \eqref{eq:linear}
\end{align} \label{eq:ots_mip}
\end{subequations}

Although model \eqref{eq:ots_mip} can be solved using off-the-shelf mixed-integer optimization solvers, such as Gurobi \cite{gurobi}, the choice of the bounds $M_{nm}, M_{mn}$ is of utmost importance. If these bounds are too loose, the relaxations performed throughout the branch-and-bound or branch-and-cut algorithms are too poor, and the total computational burden is expected to increase significantly. In all existing works that reformulate the DC-OTS problem as a mixed-integer program, these large enough constants are assumed to be symmetric, i.e., $M_{nm}=M_{mn}$. The review paper \cite{numan2023role} collects in Table 1 a summary of all proposed symmetric big-M values used in the technical literature. In particular, the authors of \cite{moulin2010transmission} propose a method to compute big-M values based on the shortest and longest paths between two nodes. This methodology has been recently revisited in \cite{fattahi2019bound}, where the authors argue that the lowest possible value of these bounds denoted by $M^{\rm OPT}_{nm}$ and $M^{\rm OPT}_{mn}$ can be obtained by solving the following bounding problems 
\begin{subequations}\label{eq:max}
\begin{align}
& M^{\rm OPT}_{nm}:=b_{nm} \times\underset{\eqref{eq:OTS_NP_Flow_S}-\eqref{eq:OTS_NP_slack}\,\cap\,\mathcal{X}_{nm}}{\max} (\theta_n-\theta_m) \\
& M^{\rm OPT}_{mn}:=b_{nm} \times\underset{\eqref{eq:OTS_NP_Flow_S}-\eqref{eq:OTS_NP_slack}\,\cap\,\mathcal{X}_{nm}}{\max} (\theta_m-\theta_n)
\end{align}
\end{subequations}
where $\mathcal{X}_{nm}:= \{\mathbf{x} \in \mathbb{B}^{\vert\mathcal{L}_{\mathcal{S}}\vert}: x_{nm} = 0\}$ imposes that the binary variable associated with the switchable line $(n,m)$ is equal to 0.
As illustrated in \cite{fattahi2019bound}, problem \eqref{eq:max} can be unbounded in power systems where switching off lines can result in isolated subnetworks. However, due to reliability and security standards, islanding in power grids is to be avoided in general and therefore, we assume that the set of switchable lines $\mathcal{L}_{\mathcal{S}}$ is such that the connectivity of the whole power network is always guaranteed by means of a spanning subgraph. The authors in~\cite{fattahi2019bound} also show that, even when $M^{\rm OPT}_{nm}$ is finite, computing it is as hard as solving the original DC-OTS problem. Therefore,  they propose an efficient methodology to find other valid bounds for \eqref{eq:linear} as follows: 

\begin{subequations}   \label{eq: big M Fattahi}  
\begin{empheq}[left={(\mathcal{SP})} \quad]{align}  
& M_{nm} = b_{nm} \hspace{-3mm} \sum_{(i,j)\in {\rm SP}_{nm}}\frac{F_{ij}}{b_{ij}}, \quad \forall (n,m) \in \mathcal{L}_{\mathcal{S}} \\ 
& M_{mn} = b_{nm} \hspace{-3mm} \sum_{(i,j)\in {\rm SP}_{mn}}\frac{F_{ij}}{b_{ij}}, \quad \forall (n,m) \in \mathcal{L}_{\mathcal{S}} 
\end{empheq} 
\end{subequations}
where ${\rm SP}_{nm}$ is the shortest path from $n$ to $m$, and ${\rm SP}_{mn}$ the shortest path from $m$ to $n$. These shortest paths are determined on a directed graph with edge costs $c_{nm}=F_{nm}/b_{nm}$ and $c_{mn}=F_{mn}/b_{nm}$ for the lines that belong to the connected spanning subgraph, and $c_{nm}=c_{mn}=\infty$ for the switchable lines. These shortest paths can be efficiently computed using Dijkstra's algorithm \cite{cormen2022introduction}. In reference \cite{fattahi2019bound}, line capacities are assumed symmetric and therefore, the big-M values computed by \eqref{eq: big M Fattahi} are also symmetric, that is, $M_{nm}=M_{mn}$. For given line capacities $\mathbf{F}$, using equations \eqref{eq: big M Fattahi} to obtain the bounds $\mathbf{M}$ for all switchable lines is denoted as $\mathbf{M} = \mathcal{SP}\left(\mathbf{F}\right)$. Among the references reviewed in \cite{numan2023role}, the methodology proposed in \cite{fattahi2019bound} is the one that leads to tighter big-M values and therefore, this approach is used here as a benchmark.  

In the next section we propose a novel methodology to compute valid bounds that are tighter than those described in~\cite{fattahi2019bound} and therefore reduce the computational burden of solving model \eqref{eq:ots_mip}. Conversely to all existing methodologies, the one we propose in this paper allows us to compute asymmetric big-M values that yield tighter mixed-integer reformulations of the DC-OTS problem. Following the assumption in~\cite{fattahi2019bound}, the proposed tighter bounds are derived under the premise that network connectivity is guaranteed by a set of non-switchable lines forming a spanning subgraph. Interested readers are referred to \cite{li2021connectivity} to delve into how the complexity of the OTS problem increases when all lines are switchable, and network connectivity is enforced through additional constraints.

\section{Bound tightening methodology} \label{sec:methodology}

\subsection{Big-M tightening}\label{sub:bigM_tightening}

The methodology proposed in this paper to find the values of the large constants $M_{nm}, M_{mn}$ is based on the following relaxations of problems \eqref{eq:max}
\begin{subequations} \label{eq:max_lin}
\begin{align}
& M_{nm}=b_{nm} \times\underset{\mathcal{R}(\mathbf{F},\mathbf{M})\,\cap\,\mathcal{X}^0_{nm}}{\max} (\theta_n-\theta_m) \\
& M_{mn}=b_{nm} \times\underset{\mathcal{R}(\mathbf{F},\mathbf{M})\,\cap\,\mathcal{X}^0_{nm}}{\max} (\theta_m-\theta_n)
\end{align}
\end{subequations}
where the feasible region defined by \eqref{eq:OTS_NP_Flow_S}-\eqref{eq:OTS_NP_slack} is replaced by the set $\mathcal{R}(\mathbf{F},\mathbf{M}):= \{( \mathbf{p}, \mathbf{\theta}, \mathbf{f}, \mathbf{x}) \in \mathbb{R}^{2\vert\mathcal{N}\vert +\vert\mathcal{L}\vert + \vert\mathcal{L}_{\mathcal{S}}\vert}:  \eqref{eq:OTS_NP_Flow_NS}-\eqref{eq:OTS_NP_slack}, \eqref{eq:linear}\}$ based on the linearization \eqref{eq:linear}. Note that the feasible region $\mathcal{R}$ depends on the parameter vectors  
$\mathbf{F}=[(F_{nm},F_{mn}), (n,m) \in \mathcal{L}]$ and $\mathbf{M}=[(M_{nm},M_{mn}), (n,m) \in \mathcal{L}_S]$. Besides, the set $\mathcal{X}^0_{nm}:= \{\mathbf{x} \in \mathbb{R}^{\vert\mathcal{L}_{\mathcal{S}}\vert}: \mathbf{0} \leq \mathbf{x} \leq \mathbf{1}, x_{nm} = 0\}$ is a relaxation of the set $\mathcal{X}_{nm}$ in which variables~$\mathbf{x}$ can take any continuous value between 0 and 1. Similar relaxed optimization problems have been used in \cite{roald2019implied,porras2022cost} to screen out redundant constraints. For valid bound values $\mathbf{M}$, it is guaranteed that $\mathcal{X}_{nm} \subset \mathcal{X}^0_{nm}$ and therefore, $M_{nm} \geq M^{\rm OPT}_{nm}$ and $M_{mn} \geq M^{\rm OPT}_{mn}$. Besides, since optimization problems in \eqref{eq:max_lin} are linear, the proposed methodology to find valid bounds for inequality constraints \eqref{eq:linear} is computationally efficient. For the remaining of this paper, we denote problems \eqref{eq:max_lin} as \textit{bounding problems} \cite{chen2015bound}. 

Although bounding problems \eqref{eq:max_lin} are easy to solve, the proposed relaxation can yield too loose bounds such that $M_{nm} \gg M^{\rm OPT}_{nm}$ and/or $M_{mn} \gg M^{\rm OPT}_{mn}$ and therefore, the computational burden of solving \eqref{eq:ots_mip} using these bounds can still be substantial. To avoid this issue, we include additional constraints to the bounding problems \eqref{eq:max_lin} so that the obtained big-M values are as tight as possible. 

In reference \cite{porras2022cost} the authors use a constraint on the generation cost of the network-constrained unit commitment problem to efficiently remove inactive constraints of the optimization model. Inspired by this idea, one may wonder  whether it is necessary to choose big-M values that guarantee the feasibility of all integer solutions, or whether it could be more effective if these bounds were tuned to also remove some feasible but suboptimal integer solutions. For the sake of intuition, let us assume that the solution of problem \eqref{eq:max} for a given switchable line indicates that the maximum angle difference is reached when the most expensive generators are producing at maximum capacity and the cheapest units are not generating anything. Most likely, the dispatch that maximizes this angle difference is much more expensive than that obtained by the the DC-OTS problem and therefore, the actual angle difference at the optimal solution of \eqref{eq:ots_mip} is probably much lower than that computed by \eqref{eq:max}. Accordingly, we define in this paper the set $\mathcal{C}:= \{\mathbf{p} \in \mathbb{R}^{\vert\mathcal{N}\vert}: \sum_n c_n p_n \leq \overline{C}\}$, where $\overline{C}$ is an upper bound on the optimal generation cost of the DC-OTS problem. The bounding problems that consider an upper bound on the production cost are then formulated as follows:

%
\begin{subequations} \label{eq:max_cost}
\begin{empheq}[left={(\mathcal{BM})} \quad]{align} 
& M_{nm}=b_{nm} \times\underset{\mathcal{R}(\mathbf{F},\mathbf{M})\,\cap\,\mathcal{X}^0_{nm}\cap\,\mathcal{C}}{\max} (\theta_n-\theta_m) \label{eq:max_cost_upper}\\
& M_{mn}=b_{nm} \times\underset{\mathcal{R}(\mathbf{F},\mathbf{M})\,\cap\,\mathcal{X}^0_{nm}\cap\,\mathcal{C}}{\max} (\theta_m-\theta_n) \label{eq:max_cost_lower}
\end{empheq}
\end{subequations}

Since the feasible regions of bounding problems \eqref{eq:max_cost} are contained in the feasible regions of \eqref{eq:max_lin}, we can guarantee that the obtained bounds are tighter than those determined in \eqref{eq:max_lin}. Besides, by setting a cap on the operational cost, the big-M values derived from \eqref{eq:max_cost} become tighter, leveraging the economic insights embedded in the maximum cost $\overline{C}$. Obviously, the tighter the value of the upper-bound cost $\overline{C}$, the smaller the feasible regions of problems \eqref{eq:max_cost}. This implies lower big-M values and the consequent reduction of the computational burden of problem \eqref{eq:ots_mip}. For given line capacities $\mathbf{F}$, big-M values $\mathbf{M}$ and upper-bound cost $\overline{C}$, using the bounding problems \eqref{eq:max_cost} to update the big-M values for all switchable lines is denoted as $\mathbf{M} = \mathcal{BM}(\mathbf{F}, \mathbf{M}, \overline{C})$.

Importantly, while the bounds computed by \eqref{eq: big M Fattahi} according to the method proposed in \cite{fattahi2019bound} are symmetric, the big-M values obtained by the bounding problems \eqref{eq:max_cost} are not symmetric in general. Another relevant point to consider is that a decrease of the big-M values associated with a specific switchable line has an impact on the feasible region $\mathcal{R}$ of the bounding problems related to the other switchable lines. As a result, it may be necessary to solve the proposed bounding problems multiple times for the entire set of switchable lines. This approach ensures that with each iteration, the big-M values will consistently decrease and lead to more refined solutions through successive iterations. 

\subsection{Line capacity tightening}\label{sub:line_tightening}

As discussed in Subsection \ref{sub:bigM_tightening}, the feasible region of problems \eqref{eq:max_cost} is reduced by imposing an upper bound on the optimal generation cost. Following this line of thought, the feasible region $\mathcal{R}$ can also be shrunk by tightening the line capacities $\mathbf{F}$. For instance, let us consider a given transmission line through which the power flow cannot exceed 100MW due to thermal limitations. However, given the location and capacity of the generating units, the demand location and variability, the network topology and parameters, and an upper bound on the generating cost, the power flow through that line is guaranteed to be always below 80MW. In such a case, we can tighten this line capacity with the following computational advantages. By reducing the capacities of the lines in the connected spanning subgraph, the big-Ms computed by \eqref{eq: big M Fattahi} also decrease. Besides, since constraint \eqref{eq:OTS_NP_Flow_limit_S} includes the product $x_{nm}F_{nm}$, reducing the capacity of switchable lines also makes model \eqref{eq:ots_mip} tighter. For these reasons, we also propose in this Section to compute the maximum feasible flows through all transmission lines as follows:

%
\begin{subequations} \label{eq:max_flow}
\begin{empheq}[left={(\mathcal{BL})} \quad]{align} 
& F_{nm}=b_{nm} \times\underset{\mathcal{R}(\mathbf{F},\mathbf{M})\,\cap\,\mathcal{X}^1_{nm}\cap\,\mathcal{C}}{\max} (\theta_n-\theta_m) \label{eq:max_cost_upper_s}\\
& F_{mn}=b_{nm} \times\underset{\mathcal{R}(\mathbf{F},\mathbf{M})\,\cap\,\mathcal{X}^1_{nm}\cap\,\mathcal{C}}{\max} (\theta_m-\theta_n) \label{eq:max_cost_lower_s}
\end{empheq}
\end{subequations}
%
\noindent where $\mathcal{X}^1_{nm}:= \{\mathbf{x} \in \mathbb{R}^{\vert\mathcal{L}_{\mathcal{S}}\vert}: \mathbf{0} \leq \mathbf{x} \leq \mathbf{1}\}$ if $(n,m) \in \mathcal{L} \setminus \mathcal{L}_S$, and $\mathcal{X}^1_{nm}:= \{\mathbf{x} \in \mathbb{R}^{\vert\mathcal{L}_{\mathcal{S}}\vert}: \mathbf{0} \leq \mathbf{x} \leq \mathbf{1}, x_{nm} = 1\}$ if $(n,m) \in \mathcal{L}_S$. Obviously, the maximum flows determined by \eqref{eq:max_flow} are always lower than or equal to the original capacities determined by thermodynamic limitations. Updating the maximum power flows through all lines in the network using the bound problems \eqref{eq:max_flow} is denoted as $\mathbf{F} = \mathcal{BL}(\mathbf{F}, \mathbf{M}, \overline{C})$. Additionally, it is worth mentioning that the reduced line capacities computed in \eqref{eq:max_flow} can also be used in \eqref{eq: big M Fattahi} to get tighter big-M values.

\begin{algorithm}
\begin{small}
\caption{Cost-driven Bound Tightening Algorithm} \label{alg:bound_tight}
\begin{algorithmic} 
\State \textbf{Input:} Original line capacities $F_{nm},F_{mn},\, \forall (n,m) \in \mathcal{L}$, demands $d_n,\, \forall n \in \mathcal{N}$, upper-bound cost $\overline{C}$, and number of iterations $K$.
\State \textbf{Initialization:} Determine SP$_{nm}$ and SP$_{mn}$ $\forall (n,m) \in \mathcal{L}_{\mathcal{S}}$. Set $k = 0$, $\mathbf{F}^0 = [(F_{nm},F_{mn}), (n,m)\in\mathcal{L}]$, and $\mathbf{M}^0 = \mathcal{SP}(\mathbf{F}^0)$. 
\begin{enumerate}[leftmargin=15pt, label={\arabic*)}]
\item Update $k \leftarrow k+1$.
\item Depending on the method, update big-M values as
\begin{enumerate}[leftmargin=15pt, label={\roman*)}]
\item $\mathbf{M}^k \leftarrow \mathcal{SP}(\mathbf{F}^{k-1})$
\item $\mathbf{M}^k \leftarrow \mathcal{BM}(\mathbf{F}^{k-1}, \mathbf{M}^{k-1}, \overline{C})$
\end{enumerate}
\item Depending on the method, update line capacities as 
\begin{enumerate}[leftmargin=15pt, label={\roman*)}]
\item $\mathbf{F}^k \leftarrow \mathbf{F}^{k-1}$
\item $\mathbf{F}^k \leftarrow \mathcal{BL}(\mathbf{F}^{k-1}, \mathbf{M}^{k-1}, \overline{C})$
\end{enumerate}
\item If $k<K$, go to step 1). Otherwise, stop.
\end{enumerate}
\State \textbf{Output:} Bounds $\mathbf{F}^K$ and $\mathbf{M}^K$.
\end{algorithmic}
\end{small}
\end{algorithm}

\subsection{Comparison procedure}\label{sub:comparison}

In summary, the method we propose in this paper starts by finding a tight upper-bound on the generating cost $\overline{C}$. Then, the big-M values and the line capacities  are iteratively tightened by solving the bounding problems \eqref{eq:max_cost} and \eqref{eq:max_flow}, respectively. Finally, the reduced bounds are used to solve the mixed-integer formulation of the DC-OTS problem \eqref{eq:ots_mip}. Algorithm \ref{alg:bound_tight} summarizes the main steps of the proposed methodology. In order to investigate the improvements derived from the proposed methodology, we compare the computational performance of the following four variations of Algorithm \ref{alg:bound_tight}:
\begin{itemize}[leftmargin=*]
    \item \textbf{S}hortest-\textbf{P}ath approach with \textbf{O}riginal line \textbf{C}apacities (SP-OC). This is the benchmark strategy proposed in \cite{fattahi2019bound} and only includes the initialization step of Algorithm \ref{alg:bound_tight}. Thus, model \eqref{eq:ots_mip} is solved with $\mathbf{F}^0$ and $\mathbf{M}^0$.
    \item \textbf{S}hortest-\textbf{P}ath approach with \textbf{R}educed line \textbf{C}apacities (SP-RC). This is an improvement of the method proposed in \cite{fattahi2019bound} that uses the reduced line capacities obtained by the bounding problems \eqref{eq:max_flow} and updates the big-M values using \eqref{eq: big M Fattahi}. Thus, this approach runs steps 3ii) and 2i) in Algorithm \ref{alg:bound_tight}, in that order, and ignores 2ii) and  3i).
    \item \textbf{B}ound \textbf{T}ightening approach with \textbf{O}riginal line \textbf{C}apacities (BT-OC). In this strategy we propose the big-M values are obtained by solving the bounding problems \eqref{eq:max_cost} with the original line capacities in all iterations. Thus, this approach runs steps 2ii) and 3i) in Algorithm \ref{alg:bound_tight}, and ignores 2i) and 3ii).
    \item \textbf{B}ound \textbf{T}ightening approach with \textbf{R}educed line \textbf{C}apacities (BT-RC). This approach reduces the big-M values and the line capacities by solving the bounding problems and then, it is expected to yield the tightest bounds. This approach runs steps 2ii) and 3ii) in Algorithm \ref{alg:bound_tight}, and ignores 2i) and 3i). 
\end{itemize}

Furthermore, in order to analyze the impact of the upper-bound cost on the proposed tightening methodology, we compare two different procedures to compute the maximum cost $\overline{C}$ to be used in the bounding problems:
\begin{itemize}[leftmargin=*]
\item \textit{Naive approach}. This approach computes an upper bound on the cost  by satisfying the total demand with the most expensive generators. By disregarding the network constraints, this cost is the solution to the linear problem:
\begin{subequations}\label{eq:cost_naive}
\begin{align}
\overline{C} = \max_{p_n}  \quad & \sum_{n} c_{n} \, p_{n} \\
\text{s.t.} \quad & \sum_n p_n = \sum_n d_n 
\end{align}
\end{subequations}
Obviously, this upper bound on the optimal generating cost does not reduce the feasibility region of the bounding problems and is just considered here for benchmarking purposes. If model \eqref{eq:cost_naive} is used to compute the upper-bound cost, the method is denoted as XX-XX-N, where XX-XX represents the bound tightening procedure SP-OC, SP-RC, BT-OC or BT-RC.

\item \textit{Heuristic approach}. The technical literature also proposes some heuristic approaches to solve the DC-OTS problem, like the greedy algorithm described in \cite{fuller2012fast}. At each step, this algorithm disconnects one switchable line at a time, computes the resulting operating cost by solving an OPF linear problem, and fixes the status of the line that leads to the lowest cost to 0. The algorithm continues until disconnecting any remaining switchable line leads to a cost increase. Although this procedure does not lead to the optimal solution of the DC-OTS, its objective function can be close enough depending on each particular case. If this heuristic approach is used to compute the upper-bound cost, the  method is denoted as XX-XX-H.
\end{itemize}

Naturally, there are alternative strategies to compute an upper bound on the operating cost other than the two discussed above. For instance, the cost that is obtained by solving a DC-OPF problem assuming that all switchable lines are connected is also a valid bound. However, this cost lies in between those computed by way of the naive and heuristic approaches. Consequently, in the numerical experiments presented in Section \ref{sec:case_study}, we only test these two. Besides, in the procedure proposed in this paper, we utilize a maximum cost constraint $\overline{C}$ as we aim to determine the optimal topology for minimizing the operational cost. Nevertheless, this methodology could be adapted to address problems where line switching decisions aim to optimize alternative objective functions, such as mitigating the risk of wildfires \cite{rhodes2020balancing, kody2022sharing}.

In the next section, we compare the performance of the four strategies described above for each upper-bound cost using different metrics. For instance, for each switchable line we compute the big-M range relative to that determined in \cite{fattahi2019bound} as follows:
\begin{equation}
    \delta^M_{nm} = 100\frac{M_{nm} + M_{mn}}{2M^0_{nm}}
\end{equation}
where $M^0_{nm}$ are the big-M values computed by \eqref{eq: big M Fattahi} with the original capacities. For instance, in the SP-OC benchmark approach, $\delta^M_{nm} =100\%$ for all switchable lines. In the remaining methods, $\delta^M_{nm} =80\%$ means that the big-M range has been reduced a 20\% for that particular switchable line. We can also compute an average value over all switchable lines as 
\begin{equation}
    \Delta^M = \frac{\sum_{nm \in \mathcal{L}_S} \delta^M_{nm}}{|\mathcal{L}_S|}
\end{equation}

Similarly, we can define the relative range of the power flows through any transmission line as:
\begin{equation}
\delta^L_{nm} = 100\frac{F_{nm} + F_{mn}}{2F^0_{nm}}    
\end{equation}
where $F^0_{nm}$ is the original line capacity. The average value is computed as
\begin{equation}
    \Delta^L = \frac{\sum_{nm \in \mathcal{L}} \delta^L_{nm}}{|\mathcal{L}|}
\end{equation}

Apart from parameters $\Delta^M$ and $\Delta^L$, we also compare the four strategies in terms of the computational burden required to solve model \eqref{eq:ots_mip} using the bounds obtained by Algorithm \ref{alg:bound_tight}. Notice that the feasible region defined by the tighter bounds is guaranteed to include the optimal solution and therefore, the optimal decisions and objective function are the same for all the methods that yield a solution in less than one hour.

\section{Case study}\label{sec:case_study}

This section provides an overview of the computational findings obtained from the various methodologies discussed in Section \ref{sec:methodology} when applied to a practical network. Our focus is on comparing the different approaches using a realistic 118-bus network, which consists of 186 lines as documented in \cite{blumsack2006network}. This network's scale is significant enough to pose challenges for current algorithms, yet it remains manageable in terms of computational complexity. Moreover, the 118-bus system is the paradigmatic network commonly employed in numerous studies on optimal transmission switching in the technical literature \cite{fisher2008optimal, kocuk2016cycle, fattahi2019bound, fuller2012fast, crozier2022feasible, hinneck2022optimal, johnson2020knearest, yang2019line, dey2022node}. Since the choice of the spanning subgraph can notably impact the computational load of the resulting OTS problem, we evaluate the performance of the proposed methodology for 100 instances with different subsets of switchable lines and demand levels $d_n$. More specifically, the spanning subgraph of 117 non-switchable lines is randomly chosen and therefore, the 69 switchable lines are different for each instance. Likewise, the nodal demand is randomly sampled using independent uniform distributions in the range $[0.9\widehat{d}_n, 1.1\widehat{d}_n]$, where $\widehat{d}_n$ is the baseline demand. In this way, we guarantee that the results reported in what follows cover a range of OTS instances of different complexity. All optimization problems have been solved using GUROBI 9.1.2 \cite{gurobi} on a Linux-based server with CPUs clocking at 2.6 GHz, 1 thread and 8 GB of RAM. In all cases, the optimality gap has been set to $0.01\%$ and the time limit to 1 hour. We used the gurobipy API for Python to solve the problems, with all Gurobi options set to default values.

Before presenting the computational results of this case study we must clarify an implementation detail of Algorithm \ref{alg:bound_tight}. In steps 2) and 3) of this algorithm the corresponding bounding problems can be solved in parallel, then reducing the final computational burden of the bound tightening procedure. However, the results of this section are obtained by solving all bounding problems \textit{sequentially} and adding up the required time to solve each linear optimization model. By doing so, the comparison of the computational burden of the different models is more informative, especially if the number of iterations is high. Besides, this strategy allows a dynamic update of the line capacities and the big-M values. That is, the bounding problem corresponding to the k-th line can be solved using the updated line capacities of the (k-1)-th line, and so on. 

After that clarification, we start this numerical analysis by fixing the number of iterations (indicated in parenthesis for each method) to 1 and comparing the four methodologies described in Section \ref{sec:methodology} combined with the two strategies to compute the upper-bound cost. Table \ref{tab:results1} collects, for each approach, the big-M and line capacity relative ranges and the computational time averaged over the 100 random instances considered including the time of solving all the bounding problems sequentially ($T^{bnd}$) and the time of solving the resulting mixed-integer DC-OTS problem ($T^{ots}$). This table also includes the number of instances that are not solved to optimality in less than one hour ($\#U$), and the maximum optimality gap among those unsolved instances. This optimality gap is provided by the solver and computed as the relative difference between the best known upper bound and the best known lower bound on the optimal objective value of the mixed-integer optimization problem.

\begin{table}[]
\renewcommand{\arraystretch}{1.5}
\centering
\begin{tabular}{lccccccc}
\hline
Method & $\Delta^M$ & $\Delta^L$ & $T^{bnd}$ & $T^{ots}$ & $\#U$ & Max gap \\
\hline
SP-OC & 100\% & 100\% & - & 672s & 14 & 0.69\% \\
\hline
BT-OC-N(1) & 68\% & 100\% & 0s & 589s & 12 & 1.69\% \\
SP-RC-N(1) & 86\% & 77\% & 1s & 480s & 9 & 0.69\% \\
BT-RC-N(1) & 68\% & 74\% & 1s & 486s & 8 & 1.04\% \\
\hline
BT-OC-H(1) & 64\% & 100\% & 2s & 298s & 5 & 0.20\% \\
SP-RC-H(1) & 83\% & 73\% & 3s & 234s & 2 & 0.47\% \\
BT-RC-H(1) & 64\% & 68\% & 3s & 167s & 1 & 0.03\% \\
\hline
\end{tabular}
\vspace{2mm}
\caption{Impact of upper-bound cost on bound tightening performance}
\label{tab:results1}
\end{table}

If we compare the benchmark SP-OC with BT-OC-N(1), we observe that the big-M values are significantly tightened, the average time is reduced by 12\%, and the number of unsolved problems is also lower. On the contrary, the maximum gap increases from 0.69\% to 1.69\%. If the line capacities are reduced using the bounding problems, approaches SP-RC-N(1) and BT-RC-N(1) also involve computational savings that amount to 28\% approximately. In any case, it seems that using the naive upper-bound cost leads to quite modest computational improvements. 

In order to improve the performance of the bound tightening proposed in this paper, the last three rows of Table \ref{tab:results1} provide the results if the upper-bound cost is obtained by the heuristic approach described in \cite{fuller2012fast}. For the 100 instances of this case study, the average and maximum error incurred by this greedy approach amounts to 2\% and 11\%, respectively. Consequently, despite the valuable insights into generating costs provided by this heuristic procedure, the obtained solutions still exhibit a significant degree of suboptimality. By analyzing these results, we realize that using a tighter upper-bound on the operating cost has a more notable impact on the computational burden of the DC-OTS problem. For instance, even if the original line capacities are considered, the approach BT-OC-H(1) is able to halve the computational time and the number of unsolved instances yielded by SP-OC. Besides, even with one iteration, the approach BT-RC-H(1) strengthens both the line capacities and big-M values by solving the proposed bounding problems and consequently reduces the computational time by 75\% and only reports one unsolved instance. 

We continue this case study by analyzing the impact of the number of iterations through a comparison of the results collected in Table \ref{tab:results2}. Naturally, increasing the iterations leads to tighter bounds. However, the bound values seem to stabilize after three or four iterations. It is also relevant to highlight that, for the same number of iterations, BT-RC always outperforms SP-RC, which clearly indicates that the big-M values obtained by the bounding problems \eqref{eq:max_cost} are tighter than those computed by \eqref{eq: big M Fattahi}, even if the line capacities are adjusted to more realistic values. In fact, the approach BT-RC-H with one iteration yields better computational results than SP-RC-H with four iterations. It is also worth mentioning that although SP-RC-H(2) has tighter bounds than SP-RC-H(1), the former involves longer computational times and a higher number of unsolved problems. This counterintuitive result is attributed to the presolving and heuristic routines integrated into the optimization solver. Finally, it is shown that the best results are provided by BT-RC-H(4), an approach that achieves a time reduction of 88\% and is able to solve all instances in less than one hour. 

\begin{table}[]
\renewcommand{\arraystretch}{1.5}
\centering
\begin{tabular}{lcccccc}
\hline
Method & $\Delta^M$ & $\Delta^L$ & $T^{bnd}$ & $T^{ots}$ & $\#U$ & Max gap \\
\hline
SP-OC & 100\% & 100\% & - & 672s & 14 & 0.69\% \\
\hline
SP-RC-H(1) & 83\% & 73\% & 3s & 234s & 2 & 0.47\% \\
SP-RC-H(2) & 77\% & 69\% & 4s & 337s & 5 & 0.37\% \\
SP-RC-H(3) & 76\% & 68\% & 5s & 266s & 4 & 0.54\% \\
SP-RC-H(4) & 75\% & 68\% & 6s & 206s & 2 & 0.25\% \\
\hline
BT-RC-H(1) & 64\% & 68\% & 3s & 167s & 1 & 0.03\% \\
BT-RC-H(2) & 56\% & 63\% & 5s & 82s & 0 & - \\
BT-RC-H(3) & 52\% & 61\% & 7s & 74s & 0 & - \\
BT-RC-H(4) & 50\% & 59\% & 8s & 72s & 0 & - \\
\hline
\end{tabular}
\vspace{2mm}
\caption{Impact of iterations on bound tightening performance}
\label{tab:results2}
\end{table}

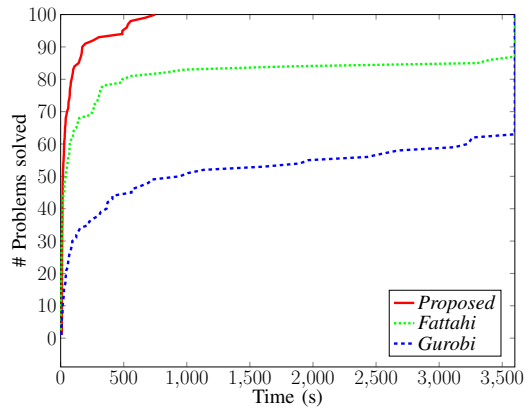
\begin{figure}
\centering
\begin{tikzpicture}[scale=0.45,font=\LARGE]
	\begin{axis}[	
    width=15cm,
    height=12cm,
    xmin = 0,
    xmax = 3600,
    ymax = 100,
    legend pos = south east,
	legend style={legend cell align=left},	
	clip marker paths=true,	
	xlabel = Time (s),
	ylabel = \# Problems solved]	
    \addplot[line width=2pt,draw=red] table [x=proposed, y=x, col sep=comma] {results.csv}; \addlegendentry{\textit{Proposed}}
    \addplot[line width=2pt,draw=green,dotted] table [x=fattahi, y=x, col sep=comma] {results.csv}; \addlegendentry{\textit{Fattahi}}
    \addplot[line width=2pt,draw=blue,dashed] table [x=gurobi, y=x, col sep=comma] {results.csv}; \addlegendentry{\textit{Gurobi}}    
	\end{axis}	
\end{tikzpicture} 
\caption{Comparison of the proposed methodology with existing benchmarks}
\label{fig:final_results}
\end{figure}

To conclude this case study, Figure \ref{fig:final_results} illustrates the number of instances solved as a function of the computational time for the following three approaches: 
\begin{itemize}
\item[-] \textit{Proposed}: This is the cost-driven bound tightening approach proposed in this paper. Among all investigated methods, we choose BT-RC-H(4) since it is the one that presents the best performance in the previous analysis.  
\item[-] \textit{Fattahi}: This methodology proposed by Fattahi et al.\ in \cite{fattahi2019bound} is based on determining the shortest-path through the spanning subgraph that connects the two nodes of every switchable line. This is the benchmark approach SP-OC that represents the state-of-the-art in the technical literature.
\item[-] \textit{Gurobi}: This strategy consists in solving the nonlinear OTS problem \eqref{eq:OTS_NP} directly with Gurobi. Gurobi is able to handle the product of binary and continuous products by using a big-M linearization with bounds that are internally computed by the solver or by adding SOS1 variables.
\end{itemize}
This figure allows us to draw the following conclusions. First, the 100 randomized instances with varying spanning subgraphs and demand levels exhibit a broad spectrum of computational difficulty. It is notable that Gurobi solves 40 ``easy'' instances in less than 500 seconds, while another 20 instances categorized as ``medium difficulty'' require computational times ranging from 500 seconds to one hour. However, 40 instances classified as ``difficult'' are not solved to optimality within one hour by this solver. Second, that the general-purpose procedure to linearize the product of binary and continuous variables implemented in Gurobi can be improved by using specific knowledge about the problem to be solved. For instance, using the power flow equations involved in the DC-OTS problem and graph theory, the shortest-path approach proposed in \cite{fattahi2019bound} provide tighter bounds than those determined internally by Gurobi to linearize the product of binary and continuous variables. Third, although the benchmark SP-OC outperforms Gurobi, the obtained big-M values can still be loose and therefore, the computational time can still be substantial for some instances. Finally, that the proposed cost-driven bound tightening methodology remarkably improves existing approaches and is able to solve the 100 random DC-OTS instances in less than 800 seconds and to reduce the average computational time by 88\% with respect to the state-of-the-art methodology.

Finally, Table \ref{tab:results3} summarizes the computational results of the proposed methodology and existing benchmarks. These results demonstrate the superiority of the bound tightening procedure we propose, with an average computational time of 80 seconds, representing a speedup of 22x and 8x compared to Gurobi or the approach described in \cite{fattahi2019bound}, respectively.

\begin{table}[]
\renewcommand{\arraystretch}{1.5}
\centering
\begin{tabular}{lccc}
\hline
Method & Time & $\#U$ & Max gap \\
\hline
\textit{Gurobi}  & 1739s & 38 & 2.03\% \\
\textit{Fattahi} & 672s & 14 & 0.69\%\\
\textit{Proposed}& 80s & 0 & -\\
\hline
\end{tabular}
\vspace{2mm}
\caption{Summary of the computational results of the proposed methodology.}
\label{tab:results3}
\end{table}



\section{Conclusions} \label{sec:conclusions}

The optimal transmission switching (OTS) aims at determining the network topology that minimizes the generating cost to satisfy a given demand. The OTS has the potential to generate substantial cost savings, but its computational requirements are high due to its typical formulation as a mixed-integer linear problem that belongs to the NP-hard class. In particular, the MIP formulation of the OTS includes big-M constants that can lead to poor relaxations if their values are too large. In this paper we propose an iterative tightening methodology that effectively reduces the big-M values, thereby alleviating the computational burden associated with the OTS. The proposed approach requires the solution of inexpensive bounding problems that account for economic information about the operating cost. Furthermore, the big-M values can be further tightened by incorporating reduced capacities of the transmission lines, which are in turn obtained by solving similar bounding problems. Using the 118-bus test system, we demonstrate that our methodology outperforms existing approaches to find proper big-M values and is able to reduce the computational of the OTS problem by 88\% in average with respect to them. 

The proposed approach is dependent on the existence of a spanning subgraph of connected lines. Therefore, extending this approach to handle the general case, where any line can be disconnected, represents a promising direction for future research.

\bibliographystyle{IEEEtran}
\bibliography{references}


  

\end{document}